# Estimating of the number of integer (natural) solutions of inhomogeneous algebraic Diophantine diagonal equations with integer coefficients

VICTOR VOLFSON

ABSTRACT. This paper investigates the upper bound of the number of integer (natural) solutions of inhomogeneous algebraic Diophantine diagonal equations with integer coefficients without a free member via the circle method of Hardy and Littlewood. Author found the upper bound of the number of natural solutions of inhomogeneous algebraic Diophantine diagonal equations with explicit variable. It was developed a method in the paper, which allows you to perform the low estimate of the number of natural (integer) solutions of algebraic Diophantine equation with integer coefficients. Author obtained a lower estimate (with this method) of the number of integer (natural) solutions for certain kinds of inhomogeneous algebraic Diophantine diagonal equations with integer coefficients with any number of variables (including Thue's equation).

1. INTRODUCTION

In general, the algebraic Diophantine equation of n-th degree with k variables and integer coefficients is $\sum_{i=0}^{n} F_i(x_1,...x_k) = 0$, where $F_i(x_1,...x_k)$ is the form of i-th degree with k variables and $F_0$ is an integer.

Thue showed that an irreducible equation $F_n(x_1,...x_k) + F_0 = 0$ (with degree $n \geq 3$ and values: $k = 2$, $F_0 \neq 0$) has only a finite number of integer solutions. Thue's method is based on approximation of algebraic numbers by rational.

Thue's equation ($n = 2$, $k = 2$) can have an infinite number of integer solutions (for example Pell's equation: $x_1^2 - Dx_2^2 = 1$, if D is not a perfect square). Thue's equation ($k \geq 2$, $n = 1$) has infinitely many solutions if there is greatest common divisor of its coefficients and free term.

---

Thue's method was developed in the works of Siegel [1], who established his famous theorem (based on this method) about a finite number of integral points on a curve of the genus $g \geq 1$. This method was extended by Schmidt [2] to the case $k > 2$, which allowed him to receive a multi-dimensional generalization of Thue's results of a finite number of integer solutions for the normative Diophantine equation.

Many methods of the theory of Diophantine equations, including Thue's method, are not effective, thus do not allow a quantitative analysis of the solutions of Diophantine equations, therefore it was developed an effective method based on the use of lower grades from the logarithms of algebraic numbers.

Currently, quantitative estimates for the integer solutions of some classical Diophantine equations obtained by this method. Baker [3] made an effective assessment of the Tue equation. A similar effective analysis for more general equation Tue-Mahler: $f(x, y) = m p_1^{x_1} ... p_s^{x_s}$ (where $p_1, ..., p_s$ are fixed numbers) was carried out in [4]. Another class of Diophantine equations (that allow effective analysis) constitutes superelliptic equation: $y^s = f(x)$, where $s \geq 2$, $f$ is an integer function of degree $n \geq 3$. Effective analysis of this equation was conducted by Baker [5]. This result was significantly strengthened [6]. Thus, the effective analysis is done only for algebraic Diophantine equations of two variables

The circular method of Hardy and Littlewood (CM) [7] allows us to make an upper estimate of the number of natural solutions to various algebraic Diophantine diagonal equations with integer coefficients and large number of variables. The method has been significantly strengthened Vinogradov [8]. It was obtained [9] an upper estimate of the number of natural solutions (with the help of the CM) of the diagonal Thue's equation: $a_1 x_1^n + ... + a_k x_k^n + a_0 = 0$ (in the case when all coefficients $a_i$ are integers and have the same sign, and $a_0$ is integer (not equal to 0) of another sign.

CM used in [10] to obtain upper estimates of the number of natural solutions of various homogeneous algebraic Diophantine diagonal equations with integer coefficients ($a_0 = 0$). It is interesting to find estimates of the number of integers (natural) for Thue's equation solutions ($k > 2$, $n \geq 2$), when the coefficients of Thue's equation have different sign.

It is also interesting to obtain lower bounds for the number of integer (natural) solutions of various inhomogeneous Diophantine diagonal equations, as a lower bound allows (in addition to the quantitative assessment) to answer the question about the existence or an infinite number of integer (natural) solution of the Diophantine equation.

2. ESTIMATING OF THE NUMBER OF INTEGERS (NATURAL) SOLUTIONS OF THE SIMPLEST INHOMOGENEOUS ALGEBRAIC DIOPHANTINE DIAGONAL EQUATIONS WITHOUT CONSTANT TERM

We make an upper bound of the number of natural solutions of the simplest inhomogeneous Diophantine diagonal equations without constant term with the help of the CM. Let us begin consideration from the homogeneous equation:

$$x_1^k - y_1^k + x_2^k - y_2^k + \ldots + x_{2^{j-1}}^k - y_{2^{j-1}}^k = 0, \qquad (2.1)$$

where $1 \leq x_i, y_i \leq N$.

It is well known Hua's lemma [7] - Suppose values: $1 \leq j \leq k$, then $\int_0^1 |f(x)|^{2^j} dx \ll N^{2^j - j + \epsilon}$, where $f(x) = \sum_{m=1}^N e^{2\pi x m^k}$. Number of natural solutions of the equation (2.1) defined by the integral [9]:

$$R_{2^j}^+(N) = \int_0^1 |f(x)|^{2^j} dx. \qquad (2.2)$$

Therefore, based on Hua's lemma we have the following upper bound for the number of natural solutions of the equation (2.1):

$$R_{2^j}^+(N) \ll N^{2^j - j + \epsilon}, \qquad (2.3)$$

where $\epsilon$ is a small positive real number and $1 \leq j \leq k$. We will consider the case $j > k$ later.

Now let us consider more general homogeneous equation:

$$x_1^k - y_1^k + x_2^k - y_2^k + \ldots + x_s^k - y_s^k = 0, \qquad (2.4)$$

where values: $1 \leq x_i, y_i \leq N$.

Based on [7] the number of natural solutions of the equation (2.4) defined by the integral:

$$R_{2s}^+(N) = \int_0^1 |f(x)|^{2s} dx \ll N^{2s - k + \epsilon}. \qquad (2.5)$$

Formula (2.5) gives an upper bound for the number of natural solutions of the equation (2.4) for values $s > 2^{k-1}$. Therefore, we can estimate the number of natural solutions of the equation (2.4) for values $k < 1 + \log_2 s$ for a fixed value $s$ according to the formula (2.5).

It was obtained [9] an estimate for the number of natural solutions of the equation (2.4) for values $k \geq 1 + \log_2 s$:

$$R_{2s}^+(N) = \int_0^1 |f(x)|^{2s}\, dx << N^{\sum_{i=1}^{t-1}(2^{j_i}-(j_i+i)/2^i)+\epsilon/2^i+2^{j_t+1}-(j_t+t)/2^{t-1}+\epsilon/2^{t-1}} = N^{2s+\epsilon-\sum_{i=1}^{t-1}(j_i+i)/2^i+(j_t+t)/2^{t-1}}. \quad (2.6)$$

Now we consider the following inhomogeneous algebraic diagonal equation:

$$x_{11}^{k_1}+...+x_{s_1}^{k_1}+...+x_{l1}^{k_l}+...+x_{ls_l}^{k_l} = y_{11}^{k_1}+...+y_{s_1}^{k_1}+...+y_{l1}^{k_l}+...+y_{ls_l}^{k_l}, \quad (2.7)$$

where $k_i$ are different natural numbers.

Based on [7] the number of solutions (2.7) corresponds to the integral:

$$\int_0^1 |f_{k_1}(x)|^{2s_1} \cdot ... \cdot |f_{k_l}(x)|^{2s_l}\, dx. \quad (2.8)$$

On the basis of (2.8) and the Cauchy-Schwarz inequality we obtain:

$$\int_0^1 |f_{k_1}(x)|^{2s_1} \cdot ... \cdot |f_{k_l}(x)|^{2s_l}\, dx \leq \left(\int_0^1 |f_{k_1}(x)|^{4s_1}\right)^{1/2} \cdot ... \cdot \left(\int_0^1 |f_{k_l}|^{2s_l^l}\, dx\right)^{1/2^{l-1}}. \quad (2.9)$$

Then, it is used formula (2.5) or formula (2.6) depending on the quantities to estimate expressions in brackets (2.9).

Thus, we recommend positing the factors in equation (2.9) in order improving accuracy of the estimate:

$$s_1 \geq s_2 \geq ... \geq s_l. \quad (2.10)$$

For example, we give an estimate of the number of natural solutions of the equation:

$$x_1^5 + x_2^4 + x_3^4 + x_4^3 + x_5^3 + x_6^3 = y_1^5 + y_2^4 + y_3^4 + y_4^3 + y_5^3 + y_6^3. \quad (2.11)$$

We consider the integral to estimate the number of natural solutions of the equation (2.11):

$$R_{12}^+(N) = \int_0^1 |f_3(x)|^6 \cdot |f_4(x)|^4 \cdot |f_5(x)|^2\, dx \leq \left(\int_0^1 |f_3(x)|^{12}\, dx\right)^{1/2} \cdot \left(\int_0^1 |f_4(x)|^{16}\, dx\right)^{1/4} \cdot \left(\int_0^1 |f_5(x)|^8\, dx\right)^{1/4}. \quad (2.12)$$

Note that in (2.12), the factors listed in the order (2.10). Then, based on (2.6) we get:

$$\int_0^1 |f_3(x)|^{12} \, dx \ll N^{12+\epsilon-(3+1)/2-(2+2)/2} = N^{8+\epsilon}, \qquad (2.13)$$

if $j_1 = 3 \le k_1 = 3$,

$$\int_0^1 |f_4(x)|^{16} \, dx \ll N^{16+\epsilon-4} = N^{12+\epsilon}, \qquad (2.14)$$

if $j_1 = 4 \le k_2 = 4$,

$$R_{12}^+(N) = \int_0^1 |f_3(x)|^6 \cdot |f_4(x)|^4 \cdot |f_5(x)|^2 \, dx \ll N^{8/2+12/4+5/4+\epsilon} = N^{8+1/4+\epsilon}, \qquad (2.15)$$

if $j_1 = 3 \le k_3 = 5$.

Based on (2.12-2.15), we get:

$$R_{12}^+(N) = \int_0^1 |f_3(x)|^6 \cdot |f_4(x)|^4 \cdot |f_5(x)|^2 \, dx \ll N^{8/2+12/4+5/4+\epsilon} = N^{8+1/4+\epsilon}. \qquad (2.16)$$

If all $k_i$ are even numbers, then (except of natural solutions of the equation (2.7)) have negative integer solutions, which are equal in absolute value. Thus the upper bound for the number of integer solutions of the equation (2.7) in this case is: $R_s(N) = 2R_s^+(N)$.

Now we consider the equation:

$$a_{11}x_{11}^{k_1} + b_{11}x_{11}^{k_1} + \ldots + a_{s1}x_{s1}^{k_1} + b_{s1}x_{s1}^{k_1} + \ldots + a_{l1}x_{l1}^{k_i} + b_{l1}x_{l1}^{k_i} + \ldots + a_{ls_1}x_{ls_i}^{k_i} + b_{ls_1}x_{ls_1}^{k_i} = 0, \qquad (2.17)$$

where all $a_{ij}, b_{ij}$ assume the value 1 or -1, and $k_l$ is an odd number.

Equation (2.17) have integer solutions $x_{ij} = -y_{ij}$, if $a_{ij}, b_{ij}$ have the same signs and solutions $x_{ij} = y_{ij}$, if $a_{ij}, b_{ij}$ have different signs. Thus, the upper estimate of the number of integer solutions of equation (2.17) coincides with the upper estimate of the number of natural solutions of the equation (2.7).

Denote $D(f)$ as the set of solutions of the Diophantine equation $f(x_1,\ldots,x_s) = 0$ in area $(Z^+)^s, s \ge 1$ (direct $s$ product of the set of natural numbers). We consider $f$ is a polynomial with integer coefficients.

It is shown [10] a simple assertion that is applicable to a wide class of Diophantine equations.

Assertion 1

Let the polynomial $f(x_1,...x_r,...x_s)$ in $(Z^+)^s$ satisfies the condition:

$$f(x_1,...x_r,...x_s) = f_1(x_1,...x_r) + f_2(x_{r+1},...x_s), \qquad (2.18)$$

then:

1. If the value $f_1 \cdot f_2 > 0$ in the rest area $(Z^+)^s$, then in $(A)^s$ (the hypercube with side $N$, dimension s) is true:

$$R_s(N) = R_r(N) \cdot R_{s-r}(N), \qquad (2.19)$$

where $R_i(N)$ is the number of solutions of the equation $f(x_1,..., x_i) = 0$ in $(A)^s$.

2. If the value $f_1 \cdot f_2 < 0$ in the rest area $(Z^+)^s$, then in $(A)^s$ is true:

$$R_s(N) \geq R_r(N) \cdot R_{s-r}(N), \qquad (2.20)$$

It follows that $R_s(N) > 0$ (from (2.19) and (2.20)), if the following conditions are satisfied:

$$R_r(N) > 0, R_{s-r}(N) > 0. \qquad (2.21)$$

Corollary 1

1. If there is a solution of the equation $f_1(x_1,...x_r) = 0$ in integers (natural) numbers and there is a solution of the equation $f_2(x_{r+1},...x_s) = 0$ in integers (natural) numbers, then the equation $f(x_1,..., x_s) = 0$ has a solution in integers (natural) numbers.

2. If there is a solution of the equation $f_1(x_1,...x_r) = 0$ in integers (natural) numbers, and there are an infinite number of solutions of the equation $f_2(x_{r+1},...x_s) = 0$ in integers (natural) numbers, then the equation $f(x_1,..., x_s) = 0$ has an infinite number of solutions in integers (natural) numbers.

It was obtained [10] the following lower estimate for the number of natural solutions of a homogeneous Diophantine equation (2.1) in $A^{2^j}$:

$$R^+_{2^j}(N) = N^{2^{j-1}} + \Sigma_{i=1}^{j-1} C^{2^i}_{2^{j-1}}(N-1)^{2^{i-1}} N^{2^{j-1}-2^{i-1}} = O(N^{2^{j-1}}). \qquad (2.22)$$

Let us consider more general homogeneous equation:

$$x^k_{11} - y^k_{11} + \ldots + x^k_{12^{j_1}} - y^k_{12^{j_1}} + \ldots + x^k_{11} - y^k_{2s1} + \ldots + x^k_{2s2^{j_k}} - y^k_{12^{j_k}} = 0, \qquad (2.23)$$

where values $1 \leq x_{ij}, y_{ij} \leq N$.

We use the fact that $2s = 2^{j_1} + \ldots + 2^{j_k}$ ($j_1 > \ldots > j_k$). Therefore, based on of Assertion 1 and (2.22) we obtain the following lower bound for the number of natural solutions of the equation (2.23) in $A^{2s}$:

$$R^+_{2s}(N) \geq \prod_{i=1}^{k} R^+_{2^{j_i}}(N) = O(N^s). \qquad (2.24)$$

For example, we give a lower estimate for the number of natural solutions of the equation:

$$x^k_{11} - y^k_{11} + x^k_{21} - y^k_{21} + x^k_{22} - y^k_{22} = 0. \qquad (2.25)$$

Based on (2.22) and (2.24) we obtain the following estimate for the equation (2.25):

$$R^+_6(N) \geq R^+_2(N) R^+_4(N) = N(2N^2 - N) = O(N^3)$$

Now let us make a lower estimate for the number of natural solutions of the homogeneous equation (2.7) in $A^s$.

Based on Assertion 1 the number of natural solutions of the inhomogeneous equation is not less than the product of natural solutions of its homogeneous parts:

$$R^+_s(N) \geq \prod_i R^+_{s_i}(N), \qquad (2.26)$$

where value $s = \sum_i s_i$.

For example, we give a lower estimate for the number of natural solutions of the equation:

$$x^2_{11} - y^2_{11} + x^4_{21} - y^4_{21} + x^4_{22} - y^4_{22} = 0. \qquad (2.27)$$

Based on (2.22) and (2.26) we obtain the following estimate for the equation (2.27):

$$R^+_6(N) \geq R^+_2(N) R^+_4(N) = N(2N^2 - N) = O(N^3)$$

Let us compare this result with the upper estimate for the number of natural solutions of this equation. Based on (2.9) we obtain the following upper bound:

$$R^+_6(N) = \int_0^1 |f_2(x)|^2 |f_4(x)|^4 \, dx \leq (\int_0^1 |f_2(x)|^4 \, dx)^{1/2} (\int_0^1 |f_4(x)|^8 \, dx)^{1/2} \ll N^{3,5+\xi}.$$

However, if in the equation: $x_1^{n_1} - x_1^{n_{s+1}} + ... + x_s^{n_s} - x_s^{n_{2s}} = 0$ has not equal degrees of the neighboring members and the condition performed: $n_1 > n_{s+1}, ... n_s > n_{2s}$, then it is easy to show that this equation has only one natural solution $x_1 = 1, ... x_s = 1$.

Indeed, if value $x_1 > 1$, we get $x_1^{n_1} > x_1^{n_{s+1}}$ and $x_1^{n_1} - x_1^{n_{s+1}} + ... + x_s^{n_s} - x_s^{n_{2s}} > 0$ etc. If $x_1^{n_s} > x_1^{n_{2s}}$, then we obtain $x_1^{n_1} - x_1^{n_{s+1}} + ... + x_s^{n_s} - x_s^{n_{2s}} > 0$ QED.

3. ESTIMATING OF THE NUMBER OF NATURAL SOLUTIONS OF THE INHOMOGENEOUS ALGEBRAIC DIOPHANTINE DIAGONAL EQUATIONS WITHOUT CONSTANT TERM

First, we consider the homogeneous diagonal equation with integer coefficients:

$$a_1 x_1^k + ... + a_s x_s^k = 0, \tag{3.1}$$

where $k$ is a positive integer, and all $a_i$ are integers.

Equation (3.1) may not have any integer solutions for small values $s$, so the lower bound of the number of integer solutions of this equation is zero. However, when the number of variables $s$ much more $k$ the situation changes.

It was proved [7] the following theorem. Suppose that $k \geq 2$ and $s_0$ is as in Vinogradov's theorem and let $s \geq min(s_0, 2^k + 1)$ and $s \geq 4k^2 - 4k + 1$. Then the equation (3.1) has non-trivial solutions in integers $x_1, ... x_s$ if we assume also that if $k$ is even then not all integers $a_1, ... a_s$ have the same sign.

Note. One can also consider (for odd values $k$) that not all $a_i$ have the same sign (if it is necessary we can replace $x_i$ on $-x_i$).

Vinogradov's theorem gives the value $s_0 < 2^k + 1$ if the value $k \geq 11$, and the value $s_0 \geq 2^k + 1$ if the value $k < 11$. Therefore, the conditions must be satisfied $s \geq 2^k + 1$ and $s \geq 4k^2 - 4k + 1$ if the value $k < 11$. For example, if values $k = 2, s \geq 15$ and values $k = 3, s \geq 34$, and only when the value $k \geq 8$ it is performed the inequality $2^k + 1 > 4k^2 - k + 1$, so enough $s \geq 2^k + 1$.

Now we consider the following inhomogeneous algebraic Diophantine diagonal equation:

$$a_{11}x_{11}^{k1} + \ldots + a_{s1}x_{s1}^{k_1} + \ldots + a_{l1}x_{l1}^{k_l} + \ldots + a_{ls_l}^{k_l}x_{ls_l}^{k_l} = 0, \tag{3.2}$$

where $k_1, \ldots k_s$, $x_{ij}$ are natural numbers and $a_{ij}$ are integers, which have different signs.

Homogeneous equation: $a_{11}x_{11}^{k1} + \ldots + a_{s1}x_{s1}^{k_1} = 0$ have integer solutions if the value $s_1, \ldots s_i$ satisfy the conditions of the above theorem. If these solutions are natural, then each homogeneous equation in the cube with a side $N$ has no less than $O(N)$ natural solutions.

Based on Assertion 1 the lower limit of the number of natural solutions of the equation (3.2) is equal to:

$$R_s^+(N) \geq O(N^i), \tag{3.3}$$

where value $s = s_1 + \ldots + s_i$.

Only certain types of homogeneous Diophantine diagonal equations have natural solutions if number of variables is small. These equations were considered in [10].

It was proved the following theorem [11]. The inhomogeneous diagonal equation without constant term:

$$a_1 x_1^{n_1} + \ldots + a_k x_k^{n_k} = 0, \tag{3.4}$$

where a natural number $k \geq 2$, and $a_1, \ldots a_k$ are integers, $a_2 + \ldots + a_k \neq 0$, $n_1, \ldots n_k$ are such natural numbers that $n_1$ and $n_2, \ldots n_k$ are relatively prime, has an infinite number of solutions in integers $x_1, \ldots x_k$ and if $a_1 > 0$, $a_2 + \ldots + a_k < 0$ it has an infinite number of solutions in natural numbers, which is given by:

$$n_1 r - (n_2 \bullet \ldots \bullet n_k) s = 1, \; t = -(a_2 + \ldots + a_k) a_1^{n_1 - 1}, \; x_1 = a_1^{n_1 \ldots n_k - 1} \bullet t^r, \; x_i = a_1^{n_1 \ldots n_k / n_i} \bullet t^{sn_2 \ldots n_k / n_i}, \tag{3.5}$$

where value $i = 2, \ldots k$.

Based on (3.5) there are determined not all solutions of (3.4), so the number of these solutions is the lower estimate.

We consider the equation: $x_1 - 2x_2^2 = 0$ as illustration. Based on (3.5) this equation has the following solutions: $x_1 = 2^{2s+1}, x_2 = 2^s$ for the natural values s. Low estimate for the number of natural solutions of this equation in the square with a side - N is: $R_2^+(N) \gg \ln(N)$. Indeed, this equation has the following number of natural solutions in a square with a side - N: $R_2^+(N) = O(\sqrt{N})$.

Let values are $a_1 > 1$ and $\sum_{i=2}^{k} a_i < -1$ in the equation (3.4), then the value t is $t = -(\sum_{i=2}^{k} a_i) \cdot a_1^{n_1-1} \neq 1$. We assume further, that the inequality $x_1 > ... > x_k$ is true. We renumber variables if it is not. Then (based on (3.6)) we will perform the following lower bound for the number of natural solutions of the equation (3.4) in a hypercube with side - $N$:

$$R_k^+(N) \geq (\ln(N) - (n_1 \bullet ... \bullet n_k - 1)\ln(a_1))/\ln(t). \quad (3.6)$$

We considered separately the case: $a_1 = 1$. Let us make the upper estimate for the number of natural solutions of the following non-homogeneous diagonal equation with explicit variable in the cube with a side - $N$:

$$x_1 = a_{n2}x_2^n + ... + a_{12}x_2 + .... + a_{nk}x_k^n + ... + a_{1k}x_k, \quad (3.7)$$

where all $a_{ij}$ are natural numbers.

First, let us estimate the number of natural solutions of the following non-homogeneous diagonal equation with explicit variable in the cube with a side - $N$:

$$x_1 = x_2^n + ... + x_k^n. \quad (3.8)$$

We prove (by induction on $k$) the estimate of the number of natural solutions of the inhomogeneous equation (3.8):

$$R_k^+(N) = O(N^{(k-1)/n}), \quad (3.9)$$

We get the equation $x_1 = x_2^n$ for the value $k = 2$, for which the estimate of the number of natural solutions to hypercube with the side - N is true: $R_2^+(N) = N^{1/n}$, that corresponds to (3.9).

Let us assume that the following estimate for the number of natural solutions of the inhomogeneous equation (3.8) for the value $k = i$ is true: $R_i^+(N) = O(N^{(i-1)/n})$, that corresponds to (3.9).

Then the estimate for the number of natural solutions of the inhomogeneous equation (3.8) for value $k = i+1$ is true:

$$R_{i+1}^+(N) = R_i^+(N) \cdot O(N^{1/n}) = O(N^{(i-1)/n}) \cdot O(N^{1/n}) = O(N^{i/n}),$$

which corresponds to (3.9) QED

We consider the estimate for the number of natural solutions of the equation $x_1 = x_2^2 + x_3^2 + x_4^2$ in the hypercube with the side - $N$ by another method.

Number for natural solutions in the hypercube with side - $N$ is defined as the number of natural solutions of the inequality: $x_1 = x_2^2 + x_3^2 + x_4^2 \leq N$, which is equal to the number of dots with natural coordinates within the sphere with radius $N^{1/2}$:

$$R_4^+(N) = \pi N^{3/2}/6 + O(N) = O(N^{3/2}). \tag{3.10}$$

Equation (3.10) corresponds to (3.9) for values $n = 2, k = 4$.

Based on the inequalities $x_2^n \leq a_{n2} x_2^n + ... + a_{12} x_2$, ... $x_k^n \leq a_{n2} x_k^n + ... + a_{12} x_k$, the number of natural solutions of the equation (3.7) in a hypercube with side $N$ will not exceed the number of natural solutions of the equation (3.8).

Therefore, we have the following upper bound for the number of natural solutions of (3.7) in a hypercube with side $N$:

$$R_k^+(N) << N^{(k-1)/n}. \tag{3.11}$$

4. ESTIMATING OF THE NUMBER OF INTEGER (NATURAL) SOLUTIONS OF THUE'S EQUATION

Let us consider the diagonal Thue's equation:

$$c_1 x_1^n + ... + c_k x_k^n + c_{k+1} = 0, \tag{4.1}$$

where $c_i (1 \leq i \leq k+1)$ - integers and k is natural number.

Equation (4.1) can be written as:

$$f = f_1 + f_2 = 0, \qquad (4.2)$$

where $f_1 = c_1 x_1^n + ... + c_{k-1} x_{k-1}^n = 0$, $f_2 = c_k x_k^n + c_{k+1} = 0$ or $f_1 = c_1 x_1^n + ... + c_{k-1} x_{k-1}^n = 0$, $f_2 = c_k x_k^n + c_{k+1} = 0$.

Equation $f_1 = 0$ is a homogeneous equation, which was studied in [10]. Equation $f_2 = 0$ (in the first case) is Thue's equation with two variables, the number of its integer (natural) solutions discussed earlier. Equation $f_2 = 0$ (in the second case) has one variable and no more than one solution for odd values n and two different integer solutions for even values n, if $|c_{k+1}/c_k|$ is an n-th power of a natural number.

Assertion 2

If the equation $f_2 = 0$ has at least one solution in integers, then the equation $f = f_1 + f_2 = 0$ has a solution in integers.

Proof of Assertion 2 follows from Assertion 1 and the fact that the equation $f_1 = 0$ always has the trivial solution $(0,...0)$.

Assertion 3

If the equation $f_1 = 0$ has at least one nontrivial solution in integers $(x_{10},...,x_{l0})$ and the equation $f_2 = 0$ has at least one solution in integers, then the equation $f = f_1 + f_2 = 0$ has infinitely many solutions in integers and has at least as much as $(2N+1)/\max\{x_{10},...,x_{l0}\}$ integer solutions in the hypercube with the side $[-N, N]$.

Proof. Equation $f_1 = 0$ is a homogeneous equation, so if it has at least one nontrivial solution in integers, then it has an infinite number of solutions in integers and it has at least as much as $(2N+1)/\max\{x_{10},...,x_{l0}\}$ integer solutions in the hypercube with the side $[-N, N]$. If the equation $f_2 = 0$ has at least one solution in integers, then (based on Assertion 1) the number of integer solutions of the equation $f = f_1 + f_2 = 0$ is infinite and it has at least as much as $(2N+1)/\max\{x_{10},...,x_{l0}\}$ integer solutions in the hypercube with the side $[-N, N]$.

Now let us consider a generalization of the Thue's equation:

$$a_{11}x_{11}^{k_1} + ... + a_{1s_1}x_{1s_1}^{k_1} + ... + a_{l11}x_{l1}^{k_{l-1}} + ... + a_{ls_l}x_{ls_l}^{k_l} + a_0 = 0, \qquad (4.3)$$

where all of the coefficients $a_{ij}$ and the free term $a_0$ are integers.

Equation (4.3) can be written as:

$$f = f_1 + f_2 = 0, \qquad (4.4)$$

where the equation $f_1 = 0$ is:

$$a_{11}x_{11}^{k_1} + ... + a_{1s_1}x_{1s_1}^{k_1} + ... + a_{l-11}x_{l-11}^{k_{l-1}} + ... + a_{l-1s_{l-1}}x_{l-1s_{l-1}}^{k_{l-1}} = 0, \qquad (4.5)$$

and the equation $f_2 = 0$ is:

$$a_{l1}x_{l1}^{k_l} + ... + a_{ls_l}x_{ls_l}^{k_l} + a_0 = 0. \qquad (4.6)$$

Equation (4.5) is homogeneous, and the equation (4.6) is Thue's equation, which is discussed above (4.1). Therefore, the Assertions 4 and 5 for the equation (4.4) are similar to Assertions 2 and 3 for the equation (4.1).

Assertion 4

If the equation $f_2 = 0$ has at least one solution in integers, then the equation (4.4) also has at least one solution in integers.

The proof of this assertion follows from the existence of the trivial solution of the equation $f_1 = 0$ and Assertion 1.

Assertion 5

If the equation $f_1 = 0$ has at least one nontrivial solution in integers $(x_{10}, ..., x_{l0})$ and the equation $f_2 = 0$ has at least one solution in integers, then the equation (4.4) has an infinite number of solutions in integers and it has at least as much as $(2N+1)/\max\{x_{10}, ..., x_{l0}\}$ integer solutions in the hypercube with the side $[-N, N]$.

Доказательство. Уравнение $f_1 = 0$ является однородным уравнением, поэтому если оно имеет хотя бы одно нетривиальное решение в целых числах, то оно имеет бесконечное число решений в целых числах и в гиперкубе со стороной $[-N, N]$ имеет не менее $(2N+1)/\max\{x_{10}, ..., x_{l0}\}$ целых решений. Если уравнение $f_2 = 0$ имеет хотя бы одно решение в целых числах, то на основании утверждения 1 количество целых решений

уравнения(41) бесконечно и в гиперкубе со стороной $[-N, N]$ имеется не менее $(2N+1)/\max\{x_{10},...,x_{l0}\}$ целых решений.

Proof. The equation $f_1 = 0$ is a homogeneous equation, so if it has at least one nontrivial solution in integers, then it has an infinite number of integer solutions and it has at least as much as $(2N+1)/\max\{x_{10},...,x_{l0}\}$ integer solutions in the hypercube with the side $[-N, N]$. If the equation $f_2 = 0$ has at least one solution in integers, then based on Assertion 1 the number of integer solutions of the equation (4.4) is infinite and it has at least as much as $(2N+1)/\max\{x_{10},...,x_{l0}\}$ integer solutions in the hypercube with the side $[-N, N]$.

Thue's equation (generalized equation) for three or more variables with the value $n \geq 3$ may have an infinite number of integer solutions unlike the case of two variables.

For example, the equation:

$$x_1^3 + x_2^3 + x_3^3 - 1 = 0 \qquad (4.7)$$

has an infinite number of integer solutions. We represent this equation as: $f = f_1 + f_2 = 0$, where $f_1 = x_1^3 + x_2^3 = 0$ and $f_2 = x_3^3 - 1 = 0$. The homogeneous equation $f_1 = 0$ has the solution $x_1 = -x_2$, and the equation $f_2 = 0$ has only one integer solution: $x_3 = 1$. Therefore, based on Assertion 1, the following lower estimate for the number of integral solutions for equation (4.7) in the hypercube with the side $[-N, N]$ is:

$$R_3(N) \geq 2N + 1 . \qquad (4.8)$$

It was found [12] an upper estimate for the number of integer solutions of the equation (4.7):

$$R_3(N) \ll N^{1+\xi} \qquad (4.9)$$

where $\xi$ is a small positive real number.

Upper estimate for the number of integer solutions (4.9) corresponds with the lower estimate of the number of integer solutions (4.8) for the equation (4.1).

Now let us look at an example for a generalized diagonal Thue's equation:
$$x_1^3 + x_2^3 + x_3^2 - 2x_4^2 - 1 = 0 \ . \ (4.10)$$

Equation (4.10) also has an infinite number of integer solutions, as it can be represented as $f = f_1 + f_2 = 0$, where $f_1 = x_1^3 + x_2^3 = 0$ and $f_2 = x_3^2 - 2x_4^2 - 1 = 0$ is Pell's equation, which has an infinite number of solutions in the case: $x_{31} = 3, x_{41} = 2, x_{3n+1} = 3x_{3n} + 2x_{4n}, x_{4n+1} = 2x_{3n} + 3x_{4n}$. Equation $f_2 = 0$ has $O(\ln N)$ integer solutions in the square with the side $[-N, N]$. Since the equation $f_1 = 0$ has $2N + 1$ integer solutions in the square with the side $[-N, N]$, thus (based on Assertion 1) we have the following low estimate for the number of integer solutions of the equation (4.10) in the hypercube with the side $[-N, N]$: $R_4(N) \gg N \ln N$.

## 5. CONCLUSION AND SUGGESTIONS FOR FURTHER WORK

The next article will continue to study estimations of the number of solutions of inhomogeneous algebraic Diophantine equations with integer coefficients.

## 6. ACKNOWLEDGEMENTS

Thanks to everyone who has contributed to the discussion of this paper.